\magnification1200
\parindent0pt
\baselineskip=1.4\baselineskip
\input pictex
\input amssym.tex

\def\tw{{\hat\vartheta}_w}

\def\ts{{\tilde\vartheta}_s}
\def\td{{\tilde\vartheta}_d}

\footline{\hss{\sl Jaworski, Zieli{\'n}ski;\ \folio}}

\def\frac#1#2{{#1\over#2}}
\def\binom#1#2{{#1\choose#2}}

\centerline{\bf THE SHORTEST}
\centerline{\bf CONFIDENCE INTERVAL FOR THE WEIGHTED SUM}
\centerline{\bf OF TWO BINOMIAL PROPORTIONS}

\vskip1truecm
\centerline{\bf Stanis³aw Jaworski}
\centerline{\bf Wojciech Zieli{\'n}ski}
\bigskip
\centerline{Department of Econometrics and Statistics}
\centerline{Warsaw University of Life Sciences}
\centerline{Nowoursynowska 159, PL-02-787 Warszawa}
\centerline{e-mail: stanislaw$\_$jaworski@sggw.pl}
\centerline{e-mail: wojciech$\_$zielinski@sggw.pl}
\centerline{http: wojtek.zielinski.statystyka.info}

\vskip1.5truecm

{\bf Summary.} Interval estimation of the probability of success in a Binomial model is considered. Zieli{\'n}ski (2018) showed that the confidence interval which uses information about non-homogeneity of the sample is better than the classical one. In the following paper the shortest confidence interval for non-homogenous sample is constructed.

\bigskip

{\bf Keywords}: confidence interval, binomial distribution, binomial proportions

{\bf AMS Classification:} 62F25

\bigskip

{\bf Introduction.} Suppose that there are two suppliers of an item with unknown defectiveness $\theta_1$ and $\theta_2$, respectively. It is known that the share of supply provided by the first supplier is $w_1$, while the share of the second one equals $w_2=1-w_1$. It is of interest to estimate the overall defectiveness $\vartheta=w_1\theta_1+w_2\theta_2$ on the basis of a sample of size~$n$.

 The common approach (for example Decrouez \& Robinson (2012)) is as follows: first we estimate $\theta_1$ and $\theta_2$ and then we construct a confidence interval for $\vartheta=w_1\theta_1+w_2\theta_2$. Zieli{\'n}ski (2018) another approach applied. We want to construct a confidence interval for $\vartheta$, and we are not interested in estimating $\theta_1$ and $\theta_2$. Note that for given $\vartheta$ there are infinitely many $\theta_1$ and $\theta_2$ giving $\vartheta$. Hence averaging with respect to $\theta_1$ and $\theta_2$ is applied. It was shown that the coverage probability of obtained confidence interval is at least the nominal confidence level. In what follows the shortest confidence interval is constructed.

We confine ourselves to the case $w_1,w_2>0$ and $w_1+w_2=1$, because of its nice interpretation mentioned above.

\bigskip
{\bf Confidence interval.}
Suppose that $n_1$ trials have been conducted with probability of success $\theta_1$, and $n_2$ trials with probability of success $\theta_2$. We are interested in estimating of $\vartheta=w_1\theta_1+w_2\theta_2$ for known $0<w_1<1$ and $w_2=1-w_1$. Assume that $w_1\leq w_2$.

Let $\xi_1\sim Bin(n_1,\theta_1)$, $\xi_2\sim Bin(n_2,\theta_2)$ and consider the random variable $$\tw=w_1\frac{\xi_1}{n_1}+w_2\frac{\xi_2}{n_2}.$$

Since we are interested in estimating $\vartheta=w_1\theta_1+w_2\theta_2$ on the basis of $\tw$, we consider the new statistical model
$$\left({\cal U},\left\{{\cal P}(n_1,n_2,\vartheta), 0\leq\vartheta\leq1\right\}\right),$$
where
$${\cal U}=\left\{w_1\frac{k_1}{n_1}+w_2\frac{k_2}{n_2}: k_1\in\{0,1,\ldots,n_1\},k_2\in\{0,1,\ldots,n_2\}\right\}.$$
The family $\left\{{\cal P}(n_1,n_2,\vartheta), 0\leq\vartheta\leq1\right\}$ of distributions  is as follows. Since for a~given $\vartheta\in(0,1)$ the probability $\theta_1$ is a number from the interval $(a(\vartheta),b(\vartheta))$, where
$$a(\vartheta)=\max\left\{0,{\vartheta-w_2\over w_1}\right\},\ b(\vartheta)=\min\left\{1,{\vartheta\over w_1}\right\}\hbox{ and }L(\vartheta)=b(\vartheta)-a(\vartheta),$$
The probability of the event $\{\tw\leq u\}$ (for $u\in{\cal U}$) equals (simply apply the law of total probability and averaging with respect to $\theta_1$)
$$
P_\vartheta\{\tw\leq u\}={1\over L(\vartheta)}\int_{a(\vartheta)}^{b(\vartheta)}\sum_{i_2=0}^{n_2}P_{\theta_1}\left\{\xi_1\leq \frac{n_1}{w_1}\left(u-\frac{w_2}{n_2}i_2\right)\right\}P_{\theta_2}\left\{\xi_2=i_2\right\}d\theta_1. \eqno{(\star)}$$
In the above integral $\theta_2=(\vartheta-w_1\theta_1)/w_2$.

Note that the family $\left\{{\cal P}(n_1,n_2,\vartheta), 0\leq\vartheta\leq1\right\}$ of distributions is decreasing in $\vartheta$, i.e. for a~given $u\in{\cal U}$,
$$P_{\vartheta_1}\{\tw\leq u\}\geq P_{\vartheta_2}\{\tw\leq u\} \quad\hbox{for}\quad \vartheta_1<\vartheta_2.$$
It follows from the fact that the family of binomial distributions is decreasing in probability of a success and $P_\vartheta\{\tw\leq u\}$ is a convex combination of binomial distributions.

Let $\tw=u$ be observed. Let $\gamma_1\in[0,1-\gamma]$. The confidence interval for $\vartheta$ at confidence level $\gamma$ is $(\vartheta_L^w(u,\gamma_1),\vartheta_U^w(u,\gamma_1))$, where
$$\eqalign{
\vartheta_L^w(u,\gamma_1)&=\cases{0&for $u=0$,\cr \max\{\vartheta: P_{\vartheta}\{\tw< u\}\}={\gamma+\gamma_1}&for $u>0$,\cr}\cr
\vartheta_U^w(u,\gamma_1)&=\cases{1&for $u=1$,\cr \min\{\vartheta: P_{\vartheta}\{\tw\leq u\}\}={\gamma_1}&for $u<1$.\cr}.\cr}$$

Note that
$$P_{\theta_1}\left\{\xi_1<t\right\}=
\cases{
P_{\theta_1}\left\{\xi_1\leq t-1\right\},& if $t$ is integer,\cr
P_{\theta_1}\left\{\xi_1\leq t\right\},& elsewhere.\cr
}$$

For given $\vartheta\in[0,1]$, the expected length of the confidence interval covering $\theta$ equals
$$l^w(\vartheta,\gamma_1)=\sum_{u\in{\cal U}}\left(\vartheta_U^w(u,\gamma_1)-\vartheta_L^w(u,\gamma_1)\right)g(u;\vartheta){\bf1}_{(\vartheta_L^w(u,\gamma_1),\vartheta_U^w(u,\gamma_1))}(\theta),$$
where
$$\eqalign{
g(u;\vartheta)&=P_\vartheta\{\tw=u\}=P_\vartheta\left\{\frac{w_1}{n_1}\xi_1+\frac{w_2}{n_2}\xi_2=u\right\}\cr
&={1\over L(\vartheta)}\int_{a(\vartheta)}^{b(\vartheta)}\sum_{i_2=0}^{n_2}P_{\theta_1}\left\{\xi_1={n_1\over w_1}\left(u-w_2{i_2\over n_2}\right)\right\}P_{{\vartheta-w_1\theta_1\over w_2}}\left\{\xi_2=i_2\right\}d\theta_1.\cr
}$$
The standard c.i. is obtained by taking $\gamma_1=(1-\gamma)/2$. Now the problem is in finding $\gamma_1$ minimizing the length $l^w(u,\gamma_1)$.

Note that it may be assumed that $u\leq0.5$. If $u>0.5$ the role of success and failure should be interchanged.

\bigskip

{\bf Theorem 1.}

{\bf 1.} For $u\leq\max\left\{\frac{w_1}{n_1},\frac{w_2}{n_2}\right\}$ the shortest confidence interval is one-sided.

{\bf 2.} In other cases the shortest confidence interval is two-sided.

\bigskip

{\bf  Lemma.} Let $f,g:[0,1]\to[0,1]$ be two continuous functions such that $f(0)=g(0)=1$, $f(1)=g(1)=0$, $f(x)>g(x)$ for all $x\in(0,1)$ and both functions are strictly decreasing.
Let $\gamma\in(0,0.5)$ and let $l(\gamma_1)=f^{-1}(\gamma_1)-g^{-1}(\gamma+\gamma_1)$ for $\gamma_1\in[0,1-\gamma]$.

{\bf a.} If $f'(0)=g'(0)=0$ and $f'(1)=g'(1)=0$ then there exists $\gamma_1^*\in(0,1-\gamma)$ such that $l(\gamma_1^*)=\min\{l(\gamma_1): \gamma_1\in[0,1-\gamma]\}$.

{\bf b.} If $f'(0)=0$, $g'(0+)=-\infty$ and $f'(1)=g'(1)=0$ then $l(1-\gamma)=\min\{l(\gamma_1): \gamma_1\in[0,1-\gamma]\}$.

{\bf c.} If $f'(0)=g'(0)=0$ and $f'(1-)=-\infty$, $g'(1)=0$ then $l(0)=\min\{l(\gamma_1): \gamma_1\in[0,1-\gamma]\}$.

\bigskip

{\it Proof of Lemma.} The derivative of $l(\gamma_1)$ with respect to $\gamma_1$ equals
$$\frac{dl(\gamma_1)}{d\gamma_1}=\frac{1}{f'(f^{-1}(\gamma_1))}-\frac{1}{g'(g^{-1}(\gamma+\gamma_1))}.$$
We have
$$\cases{
\gamma_1\to0\Rightarrow f^{-1}(\gamma_1)\to f^{-1}(0)=1 \hbox{ and }g^{-1}(\gamma+\gamma_1)\to g^{-1}(\gamma)\cr
\gamma_1\to1-\gamma\Rightarrow f^{-1}(\gamma_1)\to f^{-1}(1-\gamma) \hbox{ and }g^{-1}(\gamma+\gamma_1)\to g^{-1}(1)=0\cr
}$$
\medskip
{\bf a.} It is seen that
$$\gamma_1\to0\Rightarrow \frac{dl(\gamma_1)}{d\gamma_1}\to-\infty\hbox{ and }\gamma_1\to1-\gamma\Rightarrow \frac{dl(\gamma_1)}{d\gamma_1}\to+\infty.$$
Hence there exists $\gamma_1^*\in(0,1-\gamma)$ such that $\frac{dl(\gamma_1)}{d\gamma_1}\left|_{\gamma_1^*}\right.=0$. Since $\frac{dl(\gamma_1)}{d\gamma_1}<0$ for $\gamma_1<\gamma_1^*$ and $\frac{dl(\gamma_1)}{d\gamma_1}>0$ for $\gamma_1>\gamma_1^*$ hence we obtain thesis.
\medskip
{\bf b.} It is seen that
$$\gamma_1\to0\Rightarrow\frac{dl(\gamma_1)}{d\gamma_1}<0 \hbox{ and } \gamma_1\to1-\gamma\Rightarrow \frac{dl(\gamma_1)}{d\gamma_1}<0.$$
Hence  $\frac{dl(\gamma_1)}{d\gamma_1}<0$ for all $\gamma_1\in[0,1-\gamma]$. Therefore we obtain thesis.
\medskip
{\bf c.} It is seen that
$$\gamma_1\to0\Rightarrow \frac{dl(\gamma_1)}{d\gamma_1}\to+\infty\hbox{ and }\gamma_1\to1-\gamma\Rightarrow \frac{dl(\gamma_1)}{d\gamma_1}>0.$$
Hence  $\frac{dl(\gamma_1)}{d\gamma_1}>0$ for all $\gamma_1\in[0,1-\gamma]$. Therefore we obtain thesis.
\medskip

\bigskip

{\it Proof of Theorem 1.}

For a given $u\leq0.5$ let $G_u(\vartheta)=P_{\vartheta}\{\tw< u\}$ and $F_u(\vartheta)=P_{\vartheta}\{\tw\leq u\}$. We have
$$\vartheta_L^w(u,\gamma_1)=G_u^{-1}(\gamma+\gamma_1)\quad\hbox{ and }\quad\vartheta_U^w(u,\gamma_1)=F_u^{-1}(\gamma_1).$$
The length of the confidence interval equals
$$l^w(u,\gamma_1)=\vartheta_U^w(u,\gamma_1)-\vartheta_L^w(u,\gamma_1)=F_u^{-1}(\gamma_1)-G_u^{-1}(\gamma+\gamma_1).$$
Let
$$U(\vartheta;u,\theta_1)=\sum_{i_2=0}^{n_2}P_{\theta_1}\left\{\xi_1\leq \frac{n_1}{w_1}\left(u-\frac{w_2}{n_2}i_2\right)\right\}P_{\frac{\vartheta-w_1\theta_1}{w_2}}\left\{\xi_2=i_2\right\}$$
and
$$L(\vartheta;u,\theta_1)=\sum_{i_2=0}^{n_2}P_{\theta_1}\left\{\xi_1< \frac{n_1}{w_1}\left(u-\frac{w_2}{n_2}i_2\right)\right\}P_{\frac{\vartheta-w_1\theta_1}{w_2}}\left\{\xi_2=i_2\right\}.$$
We have$$
F_u(\vartheta)=\cases{
{w_1\over\vartheta}\int_{0}^{\vartheta/w_1}U(\vartheta;u,\theta_1)d\theta_1,&$\vartheta<w_1$,\cr
\int_{0}^{1}U(\vartheta;u,\theta_1)d\theta_1,&$w_1<\vartheta<w_2$,\cr
{w_1\over 1-\vartheta}\int_{\vartheta-w_2\over w_1}^{1}U(\vartheta;u,\theta_1)d\theta_1,&$\vartheta>w_2$\cr
}$$
and
$$
G_u(\vartheta)=\cases{
{w_1\over\vartheta}\int_{0}^{\vartheta/w_1}L(\vartheta;u,\theta_1)d\theta_1,&$\vartheta<w_1$,\cr
\int_{0}^{1}L(\vartheta;u,\theta_1)d\theta_1,&$w_1<\vartheta<w_2$,\cr
{w_1\over 1-\vartheta}\int_{\vartheta-w_2\over w_1}^{1}L(\vartheta;u,\theta_1)d\theta_1,&$\vartheta>w_2$.\cr
}$$
The derivative $dF_u(\vartheta)/d\vartheta$ for  $0<\vartheta<w_1$ equals
$$\eqalign{
&\int_{0}^{1}\frac{dU(\vartheta;u,\frac{\vartheta}{w_1}\theta_1)}{d\vartheta}d\theta_1=\cr
&-\int_0^1\sum_{i_2=0}^{n_2}  \binom{n_2}{i_2}\left(\frac{\vartheta(1 -\theta_1) }{w_2}\right)^{i_2} \left(\frac{w_2-\vartheta(1-\theta_1) }{w_2}\right)^{n_2-i_2}\cr
 &\left(\frac{\left(\frac{\vartheta  \theta_1 }{w_1}\right)^{\frac{n_1}{w_1} \left(u-\frac{w_2}{n_2}i_2 \right)} \left(1-\frac{\vartheta  \theta_1 }{w_1}\right)^{n_1-\frac{n_1}{w_1} \left(u-\frac{w_2}{n_2}i_2 \right)-1}}{{Beta}\left(n_1-\frac{n_1}{w_1} \left(u-\frac{w_2}{n_2}i_2 \right),\frac{n_1}{w_1} \left(u-\frac{w_2}{n_2}i_2\right) +1\right)}+\right.\cr
 &\left.\frac{w_1 (i_2 w_2-n_2 \vartheta(1-\theta_1))B\left(n_1-\frac{n_1}{w_1} \left(u-\frac{w_2}{n_2}i_2 \right),\frac{n_1}{w_1} \left(u-\frac{w_2}{n_2}i_2 \right)+1;1-\frac{\vartheta  \theta_1 }{w_1}\right)}
 {\vartheta(1-\theta_1 ) (w_2-\vartheta(1-\theta_1) )}\right) d\theta_1.
\cr}$$
Here $B(\cdot,\cdot;\cdot)$ is the cumulative distribution function of the beta distribution.

It is seen that
$$\lim_{\vartheta\to0}dF_u(\vartheta)/d\vartheta=0.$$

\def\konfl#1#2#3#4{{}_2F_1\left[#1,#2;#3;#4\right]}
{\bf1.} Let $\frac{w_1}{n_1}=\max\left\{\frac{w_1}{n_1},\frac{w_2}{n_2}\right\}$. Then
$$L(\vartheta;u,\theta_1)=B\left(n_1-1,1;1-\theta_1\right)\left(1-\frac{\vartheta-w_1\theta_1}{w_2}\right)^{n_2}$$
and $dG_u(\vartheta)/d\vartheta$ for  $0<\vartheta<w_1$ equals
$$\eqalign{
-w_1\frac{\Gamma(n_2+1)}{\vartheta(1-\vartheta)}&\left(\left(1-\frac{\vartheta}{w_1}\right)^{n_1+1}\konfl{1}{n_1+n_2+1}{n_2+1}{\frac{w_2}{1-\vartheta}}-\right.\cr
&\hskip4em\left. \left(1-\frac{\vartheta}{w_2}\right)^{n_2}\konfl{1}{n_1+n_2+1}{n_2+1}{\frac{w_2-\vartheta}{1-\vartheta}} \right),\cr
}$$
where $$\konfl{x}{y}{z}{t}=\sum_{j=0}^\infty\frac{(x)_j(y)_j}{(z)_j}\frac{t^j}{j!}$$
is the confluent hypergeometric function.

It is seen that
$$\lim_{\vartheta\to0}dG_u(\vartheta)/d\vartheta\to-\infty.$$
Applying point {\bf b.} of Lemma we obtain that the shortest c.i. for $\vartheta$ is one-sided.

\bigskip
For $\frac{w_2}{n_2}=\max\left\{\frac{w_1}{n_1},\frac{w_2}{n_2}\right\}$ the proof is similar. In formula $(\star)$ the role of $\theta_1$ and $\theta_2$ should be interchanged.

\bigskip
{\bf2.} It is easy to check that in other cases $\lim_{\vartheta\to0}dG_u(\vartheta)/d\vartheta=0$. From point {\bf a.} of Lemma it follows that the shortest c.i. is two-sided.

\bigskip

For $u>0.5$ we have
$$\vartheta_U^w(u,\gamma_1)=1-\vartheta_L^w(1-u,(1-\gamma)-\gamma_1)\ \hbox{ and }\ \vartheta_L^w(u,\gamma_1)=1-\vartheta_U^w(1-u,(1-\gamma)-\gamma_1).\eqno{(\star\star)}$$

\bigskip

{\bf Theorem 2.} For $u=0.5$ the standard c.i. is the shortest one.

\bigskip
{\it Proof.}

From $(\star\star)$ it is seen $\vartheta_L^w(0.5,\gamma_1)=\vartheta_L^w(0.5,(1-\gamma)-\gamma_1)$. The equality is fulfilled for $\gamma_1=(1-\gamma)/2$.

\bigskip

For $n_1=20,\ n_2=30,\ w_1=0.30$ and $\gamma=0.95$ the coverage probability of the shortest c.i. is shown in Figure~1. Note that for some probabilities $\vartheta$ the coverage probability is smaller than the nominal confidence level. This is in contradiction with the definition of the confidence interval (see Neyman (1934), Cram{\'e}r (1946), Lehmann (1959), Silvey (1970)).

\midinsert
$$\beginpicture
\sevenrm
\setcoordinatesystem units <100truemm,100truecm>
\setplotarea x from 0 to 1, y from 0.93 to 1
\axis bottom
 ticks numbered from 0 to 1 by 0.10 /
\axis left
 ticks numbered from 0.94 to 1.00 by 0.01 /
\plot "Fig_bin_krotkie.tex"
\plot 0 0.95 1 0.95 /
\put {{\bf Figure 1. }{\rm Coverage probabilities of the shortest c.i.}} at 0.50 0.92
\put {\hskip2em {\rm for $n_1=20,\ n_2=30,\ w_1=0.30$}} at 0.50 0.915
\endpicture
$$
\endinsert

\bigskip
\vfill\eject

{\bf Randomized c.i.}

To avoid the above mentioned disadvantage a randomization is introduced (c.f. Zieli{\'n}ski 2017). Let
$${\cal U}=\left\{w_1\frac{k_1}{n_1}+w_2\frac{k_2}{n_2}: k_1=0,1,\ldots,n_1, k_2=0,1,\ldots,n_2\right\}.$$
For $u\in{\cal U}$ let
$$u^{-}=\max\left\{v<u:v\in{\cal U}\right\},\ u^{+}=\min\left\{w>u:w\in{\cal U}\right\}.$$
Suppose that $\tw=u$ is observed. Consider two r.v's
$$\eta_d\sim U(0,u-u^{-})\hbox{ and }\eta_s\sim U(0,u^{+}-u)$$
and let $\ts=\tw+\eta_s$ and $\td=\tw-\eta_d$. Distributions of $\ts$ and $\td$ are easy to obtain:
$$\eqalign{
P_\theta\{\ts\leq t\}&=P_\theta\{\tw\leq \lfloor t\rfloor^-\}+\frac{\lceil t\rceil-\lfloor t\rfloor^-}{\lfloor t\rfloor-\lfloor t\rfloor^-}P_\theta\{\tw=\lfloor t\rfloor\},\cr
P_\theta\{\td\leq t\}&=P_\theta\{\tw\leq \lfloor t\rfloor\}+\frac{\lceil t\rceil}{\lfloor t\rfloor^+-\lfloor t\rfloor}P_\theta\{\tw=\lfloor t\rfloor^+\},\cr}$$
where
$$\lfloor t\rfloor=\max\left\{v\leq t:v\in{\cal U}\right\},\ \lfloor t\rfloor^-=\max\left\{v<\lfloor t\rfloor:v\in{\cal U}\right\},\ \lfloor t\rfloor^+=\min\left\{v>\lfloor t\rfloor:v\in{\cal U}\right\},$$
$$\lceil t\rceil=t-\lfloor t\rfloor,\ \lfloor t\rfloor^-=-1\hbox{ for }\lfloor t\rfloor=0,\ \lfloor t\rfloor^+=2\hbox{ for }\lfloor t\rfloor=1.$$

The shortest confidence interval $\left(\vartheta_L,\vartheta_U\right)$ at the confidence level $\gamma$ for observed $\tw=u$, $\eta_d=y$ and $\eta_s=t$ will be obtained as a solution with respect to $\vartheta$ of the following problem:
$$\cases{
\vartheta_U-\vartheta_L=\min!\cr\noalign{\vskip2pt}
P_{\vartheta_L}\left\{\tw\leq u^{-}\right\}+\frac{t}{u-u^{-}}P_{\vartheta_L}\left\{\tw= u\right\}=\gamma_2,\cr\noalign{\vskip2pt}
P_{\vartheta_U}\left\{\tw\leq u\right\}+\frac{y}{u^{+}-u}P_{\vartheta_U}\left\{\tw= u^+\right\}=\gamma_1,\cr\noalign{\vskip2pt}
\gamma_2-\gamma_1=\gamma.\cr
}$$
Note that the randomization may be simplified in the following way. Let $\eta$ be a r.v. distributed as $U(0,1)$. Than for observed $\tw=u$ and $\eta=y$ the ends of the shortest randomized c.i. are solutions of
$$\cases{
\vartheta_U-\vartheta_L=\min!\cr\noalign{\vskip2pt}
P_{\vartheta_L}\left\{\tw\leq u^{-}\right\}+yP_{\vartheta_L}\left\{\tw= u\right\}=\gamma_2,\cr\noalign{\vskip2pt}
P_{\vartheta_U}\left\{\tw\leq u\right\}+yP_{\vartheta_U}\left\{\tw= u^+\right\}=\gamma_1,\cr\noalign{\vskip2pt}
\gamma_2-\gamma_1=\gamma\cr
}$$
or equivalently
$$\cases{
\vartheta_U-\vartheta_L=\min!\cr\noalign{\vskip2pt}
P_{\vartheta_L}\left\{\tw<u\right\}+yP_{\vartheta_L}\left\{\tw= u\right\}=\gamma_2,\cr\noalign{\vskip2pt}
P_{\vartheta_U}\left\{\tw< u^+\right\}+yP_{\vartheta_U}\left\{\tw= u^+\right\}=\gamma_1,\cr\noalign{\vskip2pt}
\gamma_2-\gamma_1=\gamma.\cr
}$$
For $u=0$ we take $\vartheta_L=0$ and for $u=1$ we take $\vartheta_U=1$.

\bigskip

{\bf Theorem 3.} For $u>0$ there exists the shortest randomized c.i.

\bigskip
{\it Proof.}

The proof is similar to the proof of Theorem 1.
\bigskip

For $u>0.5$ we have
$$\vartheta_U^w(u,\gamma_1)=1-\vartheta_L^w(1-u,1-y,(1-\gamma)-\gamma_1)\ \hbox{ and }\ \vartheta_L^w(u,y,\gamma_1)=1-\vartheta_U^w(1-u,1-y,(1-\gamma)-\gamma_1)$$

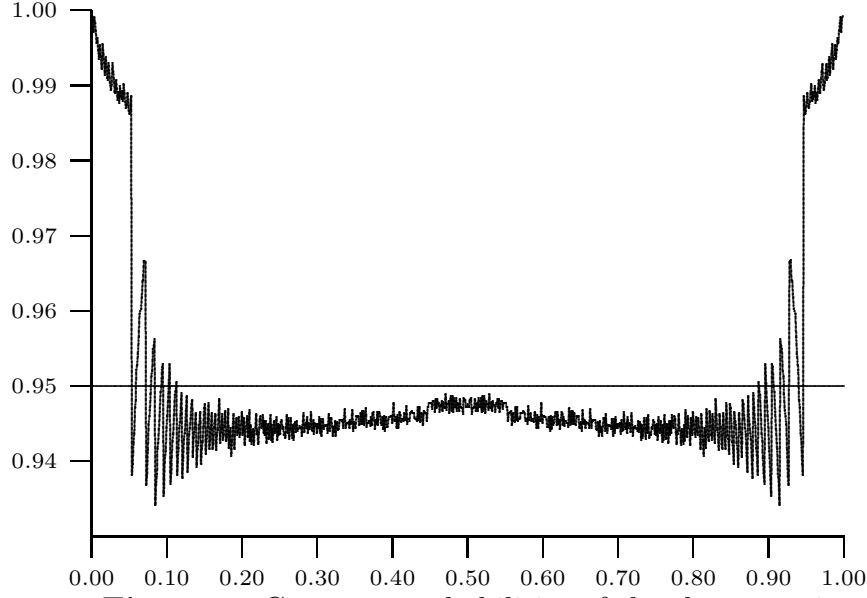
In Figure~2  the coverage probability of the shortest c.i. for $n_1=20,\ n_2=30,\ w_1=0.30$ and $\gamma=0.95$ is shown. It is seen that for all $\vartheta$ the coverage probability is at least the nominal confidence level.
\midinsert
$$\beginpicture
\sevenrm
\setcoordinatesystem units <100truemm,100truecm>
\setplotarea x from 0 to 1, y from 0.93 to 1
\axis bottom
 ticks numbered from 0 to 1 by 0.10 /
\axis left
 ticks numbered from 0.94 to 1.00 by 0.01 /
\plot "Fig_bin_krotkie_rand.tex"
\plot 0 0.95 1 0.95 /
\put {{\bf Figure 2. }{\rm Coverage probabilities of the randomized shortest c.i.}} at 0.50 0.92
\put {\hskip2em {\rm for $n_1=20,\ n_2=30,\ w_1=0.30$}} at 0.50 0.915
\endpicture
$$
\endinsert
\bigskip

{\bf An example}

Let $w_1=0.3$.
Consider an experiment consisting of $n_1=20$ and $n_2=30$ Bernoulli trials in which $\tw=0.03$ were observed. Let $\gamma=0.95$. The standard confidence interval takes on the form
$$(0.004672159,0.111730732).$$
The length of that confidence interval equals $0.107058574$.

To calculate the randomized shortest confidence interval one has to draw a value $y$ of the auxiliary variable $\eta$ and then calculate  the ends of the confidence interval. The uniform random number generator gives $y=0.2$ and the randomized shortest confidence interval takes on the form
$$(0.000883203,0.097443898).$$
The length of that confidence interval is $0.096560695$. The length of the proposed confidence interval equals $90\%$ of the length of the standard confidence interval.

The final report may look as follows:
$$w_1=0.3,\ n_1=20,\ n_2=30,\ \tw=0.03,\ y=0.2,\ \gamma=0.95,$$ $$\vartheta\in(0.000883203,0.097443898).$$

\bigskip

{\bf Conclusions}

In practical applications it is important to have conclusions as precise as possible. Hence the use of the randomized shortest confidence intervals is recommended, especially for small sample sizes. Those intervals are very easy to obtain with the aid of the standard computer software (see Appendix).


\bigskip
{\bf References}


CLOPPER, C. J. \& PEARSON, E. S. (1934), The use of confidence or fiducial limits illustrated in the case of the binomial, Biometrika 26, 404-413

CRAM{\'E}R, H. (1946). {\sl Mathematical Methods in Statistics}, Princeton University Press (Nineteenth printing 1999)

DECROUEZ, G. \&  ROBINSON, A. P. (2012) Confidence intervals for the weighted sum of two independent binomial proportions, Aust. N. Z. J. Stat. 54(3), 281-299

LEHMANN, E. L. (1959). {\sl Testing Statistical Hypothesis}, Springer (Third edition 2005)


NEYMAN, J. (1934). On the Two Different Aspects of the Representative Method: The Method of Stratified Sampling and the Method of Purposive Selection. Journal of the Royal Statistical Society 97, 558-625

SILVEY, S. D. (1970). {\sl Statistical Inference}, Chapman \& Hall (Eleventh edition 2003)



ZIELI{\'N}SKI W. (2011). Comparison of confidence intervals for fraction in finite populations, Quantitative Methods in Economy XII, 177-182.

ZIELI{\'N}SKI, W. (2016). A remark on estimating defectiveness in sampling acceptance inspection. Colloquium Biometricum 46, 9-14.

ZIELI{\'N}SKI, W. (2017) The shortest Clopper-Pearson randomized confidence interval for binomial probability. REVSTAT-Statistical Journal 15(1), 141-153.

ZIELI{\'N}SKI, W. (2018) Confidence interval for the weighted sum of two Binomial proportions. Applicationes Mathematicae 45(1), 53-60, doi: 10.4064/am2349-12-2017	

ZIELI{\'N}SKI, W. (2019) New exact confidence interval for the difference of two Binomial proportions. REVSTAT-Statistical Journal, (to appear), http://arxiv.org/abs/1903.03327.

\vfill\eject

{\bf Appendix}

An exemplary R code for calculating the shortest confidence interval is enclosed.
\bigskip

\begingroup
\obeylines
\sevenrm
\baselineskip10pt

 library(lpSolve)
  FFbin=function(m,n,q)$\scriptstyle\{$pbinom(m,n,q)$\scriptstyle\}$ \#Binomial CDF
  fbin=function(m,n,q)$\scriptstyle\{$aaa=(m==floor(m)); aaa*dbinom(floor(m),n,q)+(1-aaa)*0$\scriptstyle\}$ \#Binomial PDF
  Fbin=function(m,n,q)$\scriptstyle\{$war=(m==floor(m)); war*FFbin(m-1,n,q)+(1-war)*FFbin(m,n,q)$\scriptstyle\}$

  \#inputs
  conflevel=0.95
  w1=0.3 \#share of strata1
  n1=20 \#sample1 size
  n2=30 \#sample2 size
  k1=2 \#no of succeses in sample1
  k2=10 \#no of succeses in sample2
  \#end of inputs

  w2=1-w1
  wz=w1*k1/n1+w2*k2/n2
  los=runif(1, min = 0, max = 1)

  iksy=c(0:n2)
  smallnumber=1e-13

  \#important functions
  \#pdf at xx
  probab=function(kk,mm,qqq)$\scriptstyle\{$
   integrate(Vectorize(function(ppp)$\scriptstyle\{$fbin(kk,n1,ppp)*fbin(mm,n2,(qqq-w1*ppp)/w2)$\scriptstyle\}$),
   max(0,(qqq-w2)/w1),min(1,qqq/w1))\$value/(min(1,qqq/w1)-max(0,(qqq-w2)/w1))$\scriptstyle\}$
  density=function(xx,qqq)$\scriptstyle\{$
    gest=0
    for(ii in iksy)$\scriptstyle\{$gest=gest+probab(round((n1/w1)*(xx-(w2/n2)*ii),2),ii,qqq)$\scriptstyle\}$
    gest$\scriptstyle\}$

  \#cdf at xx
  distrib=function(xx,qqq)$\scriptstyle\{$
   integrate(Vectorize(function(ppp)$\scriptstyle\{$
   sum(FFbin(round((n1/w1)*(xx-(w2/n2)*iksy),2),n1,ppp)*fbin(iksy,n2,(qqq-w1*ppp)/w2))$\scriptstyle\}$),
   max(0,(qqq-w2)/w1),min(1,qqq/w1))\$value/(min(1,qqq/w1)-max(0,(qqq-w2)/w1))$\scriptstyle\}$

  \#randomized cdf
  randomized=function(xx,vv,yy,qqq)$\scriptstyle\{$distrib(xx,qqq)+yy*density(vv,qqq)$\scriptstyle\}$

  \#end of c.i.
  TheEnd = function(xx,vv,yy,prawd)
   $\scriptstyle\{$uniroot(function(t) $\scriptstyle\{$randomized(xx,vv,yy,t) - prawd$\scriptstyle\}$, lower = smallnumber, upper =1-smallnumber, tol = 1e-20)\$root$\scriptstyle\}$

  Leng = function(xx,ww,vv,yy,q,s)$\scriptstyle\{$TheEnd(ww,vv,yy,s)-TheEnd(xx,ww,yy,q+s)$\scriptstyle\}$

  FindMinimumLeng = function(ww,xx,vv,yy,q)$\scriptstyle\{$
  a = 0; b = 1-conflevel; \#minimum in (a,b)
  wspolczynnik = (sqrt(5) - 1)/2; \#golden ratio
  xL = b - wspolczynnik*(b - a); \#left sampling point
  xR = a + wspolczynnik*(b - a); \#right sampling point
  epsi = 10\^(-10);
  while((b- a) $\scriptstyle>$ epsi)
  $\scriptstyle\{$
    if(Leng(ww,xx,vv,los,conflevel,xL)$\scriptstyle<$Leng(ww,xx,vv,los,conflevel,xR))
    $\scriptstyle\{$b = xR; xR = xL; xL = b - wspolczynnik*(b - a);$\scriptstyle\}$
    else
    $\scriptstyle\{$a = xL; xL = xR; xR = a + wspolczynnik*(b - a);$\scriptstyle\}$
  $\scriptstyle\}$
  amin = (a + b)/2;
  amin
  $\scriptstyle\}$

  \#end of definitions

  \#determining uuplus and uuminus
  eps=0.0001
  obj.fun = c(w1/n1, w2/n2)
  constr = matrix(c(w1/n1, w2/n2, 1, 0, 0, 1) , nrow = 3, byrow = TRUE)

  constr.dir = c("$\scriptstyle<$","$\scriptstyle<$=","$\scriptstyle<$=")
  rhs = c(wz-eps,n1,n2)
  prod.sol = lp("max", obj.fun, constr, constr.dir, rhs, int.vec=1:2)
  uum=prod.sol\$objval \#objective function value

  constr.dir = c("$\scriptstyle>$","$\scriptstyle<$=","$\scriptstyle<$=")
  rhs = c(wz+eps,n1,n2)
  prod.sol = lp("min", obj.fun, constr, constr.dir, rhs, int.vec=1:2)
  uup=prod.sol\$objval \#objective function value

  uu=round(wz,8)
  uum=round(uum,8)
  uup=round(uup,8)

  smin=FindMinimumLeng(uum,uu,uup,los,conflevel)

  LowerEnd=if(uum$\scriptstyle<$=0) 0 else TheEnd(uum,uu,los,conflevel+smin)
  UpperEnd=if(uup$\scriptstyle>$=1) 1 else TheEnd(uu,uup,los,smin)

  \#output
  cat("weight of strata 1:", w1,"$\scriptstyle\backslash$n",
  "sample 1 size:", n1,"$\scriptstyle\backslash$n",
  "sample 2 size:", n2,"$\scriptstyle\backslash$n",
  "no of successes in sample 1:", k1,"$\scriptstyle\backslash$n",
  "no of successes in sample 2:", k2,"$\scriptstyle\backslash$n",
  "estimate of theta:", uu,"$\scriptstyle\backslash$n",
  "random number:", format(los,digits=3),"$\scriptstyle\backslash$n",
  "left tail:", format(smin,digits=20),"$\scriptstyle\backslash$n",
  "left end of the shortest confidence interval:", format(LowerEnd,digits=20,scientific=FALSE),"$\scriptstyle\backslash$n",
  "right end of the shortest confidence interval:", format(UpperEnd,digits=20,scientific=FALSE),"$\scriptstyle\backslash$n")

\endgroup

\bye